\documentclass[10pt,a4paper]{article}
\usepackage{amsmath,amssymb,amsfonts,amsthm}
\usepackage{float,graphicx,color}
\usepackage[all]{xy}
\usepackage{enumitem}
\DeclareMathOperator{\Div}{div}

\DeclareMathOperator{\curl}{curl}

\DeclareMathOperator{\card}{{\scriptstyle \#}}

\DeclareMathOperator{\tr}{tr}

\renewcommand{\div}{\Div}

\newcommand{\rmd}{\mathrm{d}}

\newcommand{\rmH}{\mathrm{H}}
\newcommand{\rmL}{\mathrm{L}}

\newcommand{\bbN}{\mathbb{N}}

\newcommand{\bbR}{\mathbb{R}}

\newcommand{\frh}{\mathfrak{h}}

\newcommand{\frE}{\mathfrak{E}}

\newcommand{\calH}{\mathcal{H}}

\newcommand{\calP}{\mathcal{P}}
\newcommand{\calQ}{\mathcal{Q}}

\newcommand{\calS}{\mathcal{S}}
\newcommand{\calT}{\mathcal{T}}
\newcommand{\calU}{\mathcal{U}}
\newcommand{\calV}{\mathcal{V}}

\newcommand{\from}{\leftarrow}

\newcommand{\beq}{\begin{equation}}
\newcommand{\eeq}{\end{equation}}

\newcommand{\hs}{ {\scriptscriptstyle \bigstar}}
\newcommand{\bs}{{\scriptscriptstyle \bullet}}

\newcommand{\myand}{ \quad \textrm{ and } \quad }

\newcommand{\subcell}{\lhd} % {\vartriangleleft}

\newcommand{\lspan}{\mathrm{span}}

\newcommand{\poly}{\calP}
\newcommand{\alt}{\Lambda}

\newtheorem{theorem}{Theorem}[section]
\newtheorem{lemma}[theorem]{Lemma}
\newtheorem{corollary}[theorem]{Corollary}
\newtheorem{proposition}[theorem]{Proposition}

\theoremstyle{definition}

\theoremstyle{remark}
\newtheorem{remark}{Remark}[section]

\theoremstyle{plain}

\numberwithin{equation}{section}

\begin{document}
\title{Constructions of some minimal\\ finite element systems}

\author{Snorre H.\ Christiansen\thanks{Department of Mathematics, University of Oslo, PO Box 1053 Blindern, NO-0316 Oslo, Norway. \texttt{snorrec@math.uio.no}}
~and 
Andrew Gillette\thanks{Department of Mathematics, University of Arizona, PO Box 210089, Tucson, Arizona, USA. \texttt{agillette@math.arizona.edu}}
}

%\date{}

\maketitle

\begin{abstract}
Within the framework of finite element systems, we show how spaces of differential forms may be constructed, in such a way that they are equipped with commuting interpolators and contain prescribed functions, and are minimal under these constraints. We show how various known mixed finite element spaces fulfill such a design principle, including trimmed polynomial differential forms, serendipity elements and TNT elements. We also comment on virtual element methods and provide a dimension formula for minimal compatible finite element systems containing polynomials of a given degree on hypercubes.

\end{abstract}

\section{Introduction}

A framework of Finite Element Systems (FES) has been developed, in \cite{Chr08M3AS,Chr09AWM,ChrMunOwr11,ChrRap15} to construct mixed finite elements, generalizing those of \cite{RavTho77,Ned80,BreDouMar85}. Cast in the language of differential forms, following \cite{Hip99,Chr07NM,ArnFalWin06,ArnFalWin10}, it allows for polyhedral meshes and non-polynomial differential forms.

A minimal compatible FES (mcFES) has three key properties: it is compatible, it contains certain prescribed functions, and it has the smallest dimension among all possible finite element systems with these properties.
Each property relates to a practical computational purpose: compatibility is used for the design of provably stable mixed methods, function containment is used for estimation of approximation error, and dimension minimality is used to maximize computational efficiency.

In this paper, we do the following:
\begin{itemize}
\item In Section \ref{sec:back}, we recall the main concepts of FES. We also illustrate them with some comments on the mixed Virtual Element Method \cite{BeiEtAl14}.

\item In Section \ref{sec:con}, we show how the dimension of a mcFES can be computed in terms of certain cohomology groups and how a mcFES that contains a given set of functions can be constructed, within a larger compatible FES. These results were announced, mostly in French and without proofs, in \cite{Chr10CRAS}. 

\item In Section \ref{sec:ex}, we apply this analysis and construction process to show that:
\begin{enumerate}[label=(\roman*)]
\item the trimmed polynomial spaces $\calP_r^-\Lambda^k$, defined in \cite{Hip99},  form a mcFES containing $\calP_{r-1}\Lambda^k$ on simplices;
\item the serendipity spaces $\calS_r\Lambda^k$, defined in \cite{ArnAwa14}, form a mcFES containing $\calP_{r-k}\Lambda^k$ on hypercubes;
\item the TNT elements, defined in \cite{CocQiu14}, form a mcFES containing $\calQ_r\Lambda^k$ on hypercubes;
\item the dimension of a mcFES on hypercubes that contains $\calP_r\Lambda^k$ can be given in closed form.
\end{enumerate}
\end{itemize}

\section{Background on finite element systems\label{sec:back}}

The notion of finite element systems is presented in detail in~\cite{Chr08M3AS,Chr09AWM,ChrMunOwr11,ChrRap15} so we will only recall the main definitions and key results that are most relevant to this paper.

A \emph{cellular complex} is a collection $\calT$ of cells $T$, where each cell is either a singleton or homeomorphic to the unit ball of some Euclidean space. The cells are subject to some gluing conditions that make them into a regular CW complex, but where topologists use continuous maps for these conditions, we require the maps to be at least Lipschitz continuous. We stress that in the collection $\calT$, cells of all dimensions are included. Typically, they are taken to be flat-faced polytopes, all sitting in some given $\bbR^n$, but this is not necessary for the theory to hold.
See \cite[Definition 2.1]{ChrRap15} for details.

Given a cellular complex $\calT$, a \emph{finite element system} is defined as in \cite[Definition 2.2]{ChrRap15}.
If $T$ is a cell in a cellular complex $\calT$,  we denote by $\frE^k(T)$ the set of $k$-forms on $T$ with the following property:  for any $T' \in \calT$ included in $T$ (including $T$ itself),  the pullback of the form to $T'$ is in $\rmL^2(T')$ and has its exterior derivative in $\rmL^2(T')$.

An \emph{element system} on $\calT$, is a family of  closed subspaces $E^k(T)  \subseteq \frE^k(T)$, one for each $k\in \bbN$ and each $T \in \calT$, subject to the following requirements:
\begin{itemize}
\item The exterior derivative should induce maps:
\beq
\rmd: E^k(T) \to E^{k+1}(T).
\eeq
\item If $T' \subseteq T$ are two cells in $\calT$ and  $i_{TT'}: T' \to T$ denotes the canonical injection, then pullback by $i_{TT'}$ should induce a map:
\beq
i_{TT'}^\star: E^k(T) \to E^k(T').
\eeq
\end{itemize}
For instance the spaces $\frE^k(T)$ constitute an element system. A finite element system (FES) is one in which all the spaces are finite dimensional.

We define $E^k(\calT)$ as follows :
\begin{equation}
E^k(\calT)= \{u \in \bigoplus_{T \in \calT} E^k(T) \ : \  \forall T, T' \in \calT \quad  T' \subseteq T \Rightarrow u_{T}|_{T'} = u_{T'}\}. 
\end{equation}
In this definition, which can be interpreted as encoding a continuity property of differential forms,  $u_{T}|_{T'}$ denotes the pullback of $u_T$ to $T'$ by the inclusion map.

Not all finite element systems yield good spaces $E^k(\calT)$. As in \cite[Definition 2.3]{ChrRap15}, we consider the following two conditions on an element system $E$ on a cellular complex $\calT$: 
\begin{itemize}
\item \emph{Extensions.} For each $T\in \calT$ and $k\in \bbN$, the restriction operator (pullback to the boundary) $E^k(T) \to E^k(\partial T)$ is onto. The kernel of this map is denoted $E^k_0(T)$.
\item \emph{Local exactness.} The following sequence is exact for each $T\in \calT$:
\begin{equation}\label{eq:coh}
\xymatrix{
0 \ar[r] & \bbR \ar[r] & E^0(T) \ar[r]^-\rmd & E^1(T) \ar[r]^-\rmd & \cdots \ar[r]^-\rmd &  E^{\dim T}(T) \ar[r] & 0.
}
\end{equation}
The second arrow sends an element of $ \bbR$ to the constant function on $T$ taking this value.
\end{itemize}
We will say that an element system \emph{admits extensions} if the first condition  holds, is \emph{locally exact} if the second condition holds and is \emph{compatible} if both hold.

In \cite[Proposition 2.6]{ChrRap15}, it is shown that for finite element systems admitting extensions, local exactness is equivalent to the combination of the following properties:
\begin{itemize}
\item For each $T \in \calT$,  $E^0(T)$ contains the constant functions.
\item For each $T\in \calT$, the following sequence (with boundary condition) is exact:
\begin{equation} \label{eq:coh0}
\xymatrix{
0 \ar[r] & E^0_0(T) \ar[r] & E^1_0(T) \ar[r] & \cdots \ar[r] & E^{\dim T}(T) \ar[r] & \bbR \ar[r] & 0.
}
\end{equation}
The second to last arrow is integration.
\end{itemize}

We take it for granted that a good finite element method consists of defining a compatible finite element system that contains certain prescribed functions and with which one can somehow compute.
That we want a \emph{compatible} finite element system is justified for instance in \cite[Section 2.4]{ChrRap15},  where the underlying conditions are related to the existence of degrees of freedom and commuting interpolators, the main tools of analysis of mixed methods. Typically, one would also want the finite element system to contain polynomials of a certain degree, to ensure corresponding best approximation properties in Sobolev spaces. We also remark that in some cases it is more desirable to contain certain exponentials \cite{Chr13FoCM,ChrHalSor14}.

As was already the case in \cite{Chr08M3AS}, one of the basic tools of our constructions is harmonic extension:
\begin{proposition}\label{prop:intext} Suppose $E$ is a FES where each $E^k(T)$ is equipped with a scalar product, denoted $a$. We suppose that $T$ is a cell such that  (\ref{eq:coh0}) is exact.
For each  $\alpha \in \bbR$ there is a unique $u$ in $E^{\dim T}(T)$ such that:
\begin{equation}
\int_T u = \alpha \myand \forall v \in E^{\dim T -1}_0(T)\quad  a(u, \rmd v ) = 0.
\end{equation}
Fix $k<\dim T$. Any $u \in E^k(\partial T)$ that has an extension in $E^k(T)$, has a unique extension in $ E^k(T)$ such that:
\begin{equation}\label{eq:harm}
\forall v \in E^k_0(T) \quad a(\rmd u, \rmd v)=0
\myand
\forall v \in E^{k-1}_0(T) \quad a(u, \rmd v)=0.
\end{equation}
\end{proposition}
An element $u$ of $E^k(T)$ such that (\ref{eq:harm}) holds will be called \textit{$E$-harmonic}.  The proposition above asserts that elements having an extension have a unique $E$-harmonic extension. Orthogonality with respect to $a$ will be denoted $\perp$, so that (\ref{eq:harm}) can also be written:
\begin{equation}
\rmd u \perp \rmd  E^k_0(T) \myand u \perp \rmd  E^{k-1}_0(T).
\end{equation}

We are also interested in dimension counts, so we recall \cite[Proposition 2.1]{ChrRap15}:
\newpage
\begin{proposition}\label{prop:extdim}
Let $E$ be a FES on a cellular complex $\calT$. Then: 
\begin{itemize}
\item We have:
\begin{equation}\label{eq:bounddim}
\dim E^k (\calT) \leq \sum_{ T\in \calT}\dim E^k_0(T).
\end{equation}
\item Equality holds in (\ref{eq:bounddim})  if and only if $E$ admits extensions for $k$-forms on each $T \in \calT$.
\end{itemize}
\end{proposition}

Finally we recall some results on tensor products, following \cite{Chr09AWM,ChrMunOwr11}.

Suppose $\calU$ and $\calV$ are cellular complexes, equipped with FES systems $B$ and $C$, respectively. Consider the product cellular complex:
\begin{equation}
\calT = \{ U \times V \ : \ U \in \calU \textrm{ and } V \in \calV \}.
\end{equation}
For $U \in \calU$ and $V\in \calV$, define spaces $A^k(U\times V)$ by:
\begin{equation}\label{eq:tensordef}
A^k(U\times V) = \{ p_U^\star u \wedge p_V^\star v \ : \ u \in B^l(U), \quad v \in C^{k-l}(V), \quad 0 \leq l \leq k \},
\end{equation}
where $p_U: U\times V \to U$ and $p_V : U\times V \to V$ are the canonical projections.

One checks that  $A$ is a FES on $\calT$, and we call it the tensor product of $B$ and $C$. This name is motivated by the fact that the formula:
\begin{equation}
u \otimes v = p_U^\star u \wedge p_V^\star v 
\end{equation}
indeed defines a tensor product.  Then (\ref{eq:tensordef}) can be written in compact form:
\begin{equation}
A^\bs(U\times V) = B^\bs(U) \otimes C^\bs(V).
\end{equation}
One recognizes a tensor product of graded spaces.

The important result is:
\begin{proposition}\label{prop:tensorfes}

The tensor product construction satisfies:
\begin{itemize}
\item If $B$ and $C$ admit extensions, then so does $A$.
\item If $B$ and $C$ are locally exact, then so is $A$.
\item If $B$ and $C$ are compatible, then so is $A$.
\end{itemize}
\end{proposition}
These facts were proved in \cite{Chr09AWM} and with some simplifications in \cite{ChrMunOwr11}. It is possible to simplify the proof of the first bullet point further by using the only if part of Proposition \ref{prop:extdim}.

\subsection{Comments on virtual element methods}

The Virtual Element Method is presented for scalar function spaces in \cite{BeiEtAl13}. It has recently been extended to the mixed setting, producing $\rmH(\curl)$ and $\rmH(\div)$ conforming spaces \cite{BeiEtAl14}. As we see it, the mixed virtual element method does essentially two things: it defines a finite element system and it provides a way of computing with it, that avoids reconstructing all the basis functions from the degrees of freedom. As an illustration of the methods of FES, we provide some details on the first point. The second point however, we leave open.

The constructions of \cite{Chr08M3AS} were based on a notion of harmonic extension for differential forms (both continuous and discrete). Canonical degrees of freedom for compatible FES can be defined, as in projection based interpolation, by non-homogeneous harmonic extensions -- see \cite[Proposition 3.23]{Chr09AWM} and \cite[Proposition 5.44]{ChrMunOwr11}. As we shall see, the mixed VEM spaces are defined by a similar technique.

\begin{lemma}\label{lem:zdef}
Let $T$ be a cell of dimension $n$.
\begin{itemize}
\item Fix $k <n$. Choose $f\in \rmL^2 \alt^k(T)$ and $g\in \rmL^2 \alt^{k-1} (T)$. The system, with unknown $u \in \frE^k_0(T)$:
\begin{align}
\rmd^\star \rmd u & = f,\\
\rmd^\star  u & = g.
\end{align}
has a solution if and only if $\rmd^\star f = 0$ and $\rmd^\star g = 0$. In this case the solution is unique.
\item The case $k = n$. Fix $f \in \bbR$ and $g\in \rmL^2 \alt^{k-1} (T)$. The system, with unknown $u \in \frE^k(T)$:
\begin{align}
\int u & = f,\\
\rmd^\star  u & = g,
\end{align}
has a solution if and only if $\rmd^\star g = 0$. In this case the solution is unique.
\end{itemize}
\end{lemma}
\begin{proof}
Left to the reader.
\end{proof}

Consider now a cellular complex $\calT$. We suppose we have, for each integer $k$ and cell $T \in \calT$, with $k < \dim T$, a finite dimensional subspace $Z^k(T)$ of $\rmL^2 \alt^k(T)$, whose elements $g$ satisfy $\rmd^\star g = 0$. For $k = \dim T$, we put $Z^k(T) = 0$ and note that for $u \in \rmL^2\alt^k(T)$ we have $ \rmd^\star \rmd u = 0$. We define:
\begin{equation}
A^k(T) = \{ u \in \frE^k(T) \ : \ \forall~T' \subcell T  \quad \rmd^\star \rmd u|_{T'} \in Z^k(T') \textrm{ and } \rmd^\star u|_{T'} \in Z^{k-1}(T') \},
\end{equation}  
where $T'\subcell T$ indicates that $T'$ is a face of $T$.
\begin{proposition}
The above defined spaces $A^k(T)$ constitute a compatible finite element system and
for $k < \dim T$ :
\begin{equation}
\dim A^k_0(T) = \dim Z^k(T) + \dim Z^{k-1}(T),
\end{equation}
whereas for $k = \dim T$:
\begin{equation}
\dim A^k_0(T) = 1 + \dim Z^{k-1}(T).
\end{equation}
\end{proposition}
\begin{proof}
The definition of $A^k(T)$ ensures that it is an element system.
The technique of harmonic extension shows that $A$ has the extension property.
The local exactness in the form (\ref{eq:coh0}) is also trivial to check.
Thus $A$ defines a compatible element system. 
The dimension counts follow from Lemma \ref{lem:zdef}.
\end{proof}

There are several variants of mixed VEM, but the main one seems to correspond to the following choice of $Z$. We suppose that each cell is flat, so that polynomials are well defined objects. Fix $r \geq 1$ and put, for $k < \dim T$.
\begin{equation}
Z^k(T) = \{f \in \poly_{r-1} \alt^k(T) \ : \ \rmd^\star f = 0 \}.
\end{equation}
Then $A$, defined as above, is a compatible FES by the preceding proposition, and $A^k(T)$ contains polynomial $k$-forms of degree $r$.
There is a notion of degrees of freedom for FES. Since it is rather intuitive we omit the definition; details can be found in \cite[Section 2.4]{ChrRap15}.

\begin{proposition}
For any integer $k$ and cell $T$ of dimension $n$, consider the linear forms, defined on $k$-forms by:
\begin{equation}
u \mapsto \int u \wedge v,
\end{equation}
for some $v \in \poly^-_r \alt^{n-k}(T)$.

These linear forms constitute unisolvent degrees of freedom on the FES $A$.
\end{proposition}

\begin{proof}
Let $\hs$ denote the Hodge star operator.
We notice that $u \in \poly_{r-1} \alt^k(T)$ iff $\hs u \in \poly_{r-1} \alt^{n-k}(T)$. Moreover, for such $u$ we have:
\begin{equation}
\rmd^\star u = 0 \iff \rmd \hs u = 0.
\end{equation}

Therefore, for $k < \dim T = n$ (recalling \cite[Theorem 5.4]{ArnFalWin10}):
\begin{align}
\dim Z^k(T) &= \dim \{ u \in \poly_{r-1} \alt^{n-k}(T) \ : \ \rmd u = 0\},\\
& = \dim \{ u \in \poly^-_{r} \alt^{n-k}(T) \ : \ \rmd u = 0\}.
\end{align}
We also have:
\begin{align}
\dim Z^{k-1}(T) & = \dim \{ u \in \poly_{r-1} \alt^{n-k+1}(T) \ : \ \rmd u = 0\},\\
& = \dim \{ \rmd v : v \in \poly^-_{r} \alt^{n-k}(T)\}.
\end{align}
It follows that:
\begin{equation}
\dim Z^k(T) + \dim Z^{k-1}(T) = \dim \poly^-_{r} \alt^{n-k}(T).
\end{equation}
In the case $ k= \dim T$ we may also check:
\begin{equation}
1 + \dim Z^{k-1}(T) = \dim \poly^-_{r} \alt^{n-k}(T).
\end{equation}
Therefore:
\begin{equation}
\dim A^k_0(T) = \dim \poly^-_{r} \alt^{n-k}(T).
\end{equation}
Suppose now that $u \in A^k_0(T)$ satisfies:
\begin{equation}
\forall v \in \poly^-_{r} \alt^{n-k}(T) \quad \int u \wedge v = 0.
\end{equation}
We have $\rmd^\star \rmd u \in \poly_{r-1} \alt^k(T)$, hence $\rmd \hs \rmd u \in \poly_{r-1} \alt^{n-k}(T) \subseteq \poly^-_{r} \alt^{n-k}(T)$ . Hence:
\begin{equation}
0 = \int u \wedge \rmd \hs \rmd u = \pm \int \rmd u \wedge \hs \rmd u = \pm \int |\rmd u|^2.
\end{equation}
Therefore $\rmd u = 0$.

Let $\frh$ denote the standard homotopy operator used to prove the Poincar\'e lemma, as in e.g.\ \cite{Tay96I}. Recall that $\frh$, which is proportional to $\kappa$ on homogeneous polynomial differential forms, maps $\poly^-_r \alt^{n-k+1}(T)$ to $\poly^-_r \alt^{n-k}(T)$.
We have $\rmd^\star u \in \poly_{r-1} \alt^{k-1}(T)$, so that $ \rmd \hs u \in \poly_{r-1} \alt^{n-k+1}(T)$ and hence $ \frh \rmd \hs u \in \poly^-_{r} \alt^{n-k}(T)$. We write:
\begin{align}
0 = \int u \wedge  \frh \rmd \hs u & =  - \int u \wedge \rmd \frh \hs u  + \int u \wedge \hs u ,\\
 & =  \pm \int \rmd u \wedge \frh \hs u  + \int |u|^2=   \int |u|^2.
\end{align}
Hence $u = 0$.

It follows that the integrated wedge product is an invertible bilinear form on the product $A^k_0(T) \times \poly^-_r \alt^{n-k}(T)$. This concludes the proof.
\end{proof}

\section{Construction of minimal finite element systems\label{sec:con}}

The material presented in this section is an expanded version in English of material that appeared previously, in French and without proofs, in~\cite{Chr10CRAS}.
Given any sequence
\[
\xymatrix{
\cdots \ar[r]^-\rmd & X^{k-1}(T) \ar[r]^-\rmd & X^k(T) \ar[r]^-\rmd & X^{k+1}(T) \ar[r]^-\rmd & \cdots
}
\]
satisfying $\rmd \circ \rmd=0$, we use the notation $\ker \rmd | X^k$ to denote the kernel of the $\rmd$ map whose domain is $X^k$ and $\rmd X^k$ to denote the image of that map.
The $k$-th cohomology group associated to this sequence is the quotient: 
\[
\rmH^k(X^\bs) =(\ker \rmd | X^k) /\rmd X^{k-1}.
\]
\subsection{Sufficient conditions for minimality}
\begin{proposition}
Suppose that $A$ is a finite element system, and that $B$ is a compatible finite element system containing $A$. Then we have:
\begin{equation}\label{eq:mindim}
\dim B^k_0(T) \geq \dim A^k_0(T) + \dim \rmH^{k+1}(A^\bs_0(T)).
\end{equation}
\end{proposition}
\begin{proof}
By the rank-nullity theorem and the definition of cohomology groups, we have that:
\begin{align}
\dim A^k_0(T) & = \dim \rmd A^k_0(T) + \dim \ker \rmd | A^k_0(T),\\
\dim \rmH^{k+1}(A^\bs_0(T)) & = \dim \ker \rmd | A^{k+1}_0(T) - \dim \rmd A^k_0(T) .
\end{align} 
Hence:
\begin{equation}\label{eq:comp}
\dim A^k_0(T) + \dim \rmH^{k+1}(A^\bs_0(T)) = \dim \ker \rmd | A^{k+1}_0(T) + \dim \ker \rmd | A^k_0(T).
\end{equation}
Now, since $B^\bs_0(T)$ is exact we may choose a subspace $V$ of $B^k_0(T)$ such that the map:
\begin{equation}
\rmd : V \to \ker \rmd | A^{k+1}_0(T),
\end{equation}
is an isomorphism. Then it is clear that:
\begin{equation}
V \cap \ker \rmd | A^k_0(T) = 0.
\end{equation}
It follows that:
\begin{equation}
\dim B^k_0(T) \geq \dim V + \dim \ker \rmd | A^k_0(T).
\end{equation}
Since $\dim V=\dim \ker \rmd | A^{k+1}_0(T)$, we apply (\ref{eq:comp}) to complete the proof.
\end{proof}

We will use this as follows:

\begin{corollary}\label{cor:minnec}
Suppose that $A$ is a finite element system, and that $B$ is a compatible finite element system containing $A$. Suppose that:
\begin{equation}\label{eq:mindimeq}
\dim B^k_0(T) = \dim A^k_0(T) + \dim \rmH^{k+1}(A^\bs_0(T)).
\end{equation}
Then $B$ is minimal among compatible finite element systems containing $A$.
\end{corollary}

We will show in the next two sections that, given a finite element system $A$, there exists a compatible finite element system $B$ containing $A$ such that (\ref{eq:mindimeq}) holds. Thus the condition is not only sufficient, it is necessary for minimality. Moreover, in the case that $A$ is included in a compatible finite element system $B$, it is always possible to find a minimal one containing $A$, inside $B$. Finally, in this context, we will actually \emph{construct} a particular minimal FES. We will show how other authors have proposed different definitions of finite element spaces, which turn out to be minimal among those containing certain differential forms.

\subsection{First step \label{sec:consminfir}}

Given a finite element system $A$, the first step is to construct a finite element system $\tilde A$ containing $A$ such that $\tilde A_0$ is exact.
In the case that $A_0$ is exact already, this step can be skipped.

For $k < \dim T$ choose a subspace $H^{k+1}(T)$ such that:
\begin{equation}
\rmd A^k_0(T) \oplus H^{k+1}(T) =   \ker \rmd | A^{k+1}_0(T).
\end{equation}
Then choose a subspace $E^k(T)$ of $B^k_0(T)$ such that $\rmd$ restricts to an isomorphism:
\begin{equation}
\rmd : E^k(T) \to H^{k+1}(T).
\end{equation}
For $k=\dim T$, if $A^k(T)$ contains an element with integral $1$, set $E^k(T)= 0$, if not choose an element in $B^k(T)$ with integral $1$ and let $E^k(T)$ be its linear span.

We remark that $A^k(T) \cap E^k(T) = 0$ and set: 
\begin{equation}
\tilde A^k(T) := A^k(T) \oplus E^k(T).
\end{equation}
We remark that we have an induced isomorphism: 
\begin{equation}\label{eq:disoe}
\rmd : E^k(T) \to \rmH^{k+1}(A^\bs_0(T)),
\end{equation}
where if $k=\dim T$, the exterior derivative operation is replaced by integration. In fact we might take this, together with $E^k(T) \subseteq B^k_0(T)$,  as the defining properties of $E^k(T)$, and the space $\rmH^k(T)$ introduced above merely as a device to make it clear that such an $E^k(T)$ exists. 

We could also have specified $E^k(T)$ from a choice of scalar products as follows:
If $k < \dim T$, set:
\begin{equation}\label{eq:defe}
E^k(T) := \{ u \in B^k_0(T) \ : \ \rmd u \in A^{k+1}_0(T), \ \rmd u \perp \rmd A^k_0(T) \textrm{ and } u \perp \rmd B^{k-1}_0(T) \}.
\end{equation}
If $k=\dim T$ and if $A^k(T)$ contains an element with integral $1$, set $E^k(T):= 0$, if not set:
\begin{equation}\label{eq:defebis}
E^k(T):= \{ u \in B^k(T) \ : \ u \perp \rmd B^{k-1}_0(T) \}.
\end{equation}
Then one checks immediately that $E^k(T)$ is a subspace of $B^k_0(T)$  such that (\ref{eq:disoe}) is an isomorphism.

\begin{proposition}
Under the above hypothesis, $\tilde A$ is a FES containing $A$ such that the sequences:
\begin{equation}
0 \to \tilde A^0_0(T) \to \tilde A^1_0(T) \to \cdots \to \tilde A^{\dim T}(T) \to \bbR \to 0,
\end{equation}
are exact.
\end{proposition}

\begin{proof}
That $\tilde A$ is stable under restrictions and the exterior derivative is immediate, so we only need to prove the exactness property.

Consider first $k < \dim T$. 
Suppose $u \in \tilde A^{k}_0(T) $ is such that $\rmd u = 0$. We have $\tilde A^k_0(T) = A^k_0(T) \oplus E^k(T)$, so we can write $u = v + w$ with $v \in A^k_0(T)$ and $w\in E^k(T)$. Then $\rmd w = - \rmd v \in \rmd A^k_0(T)$ so $\rmd w = 0$, so $w= 0$.

Now we have $\rmd v = 0$.  If $k = 0$ this gives $v = 0$. If $k >0$,  we have  already remarked that $\rmd  A^{k-1}_0(T) \oplus \rmd E^{k-1} (T) =   \ker \rmd | A^{k}_0(T)$, so we can write $v = \rmd t$ with $t \in \tilde A^{k-1}_0(T)$. This concludes the case $k < \dim T$.

If $k = \dim T$, the same proof goes through, replacing the exterior derivative by the integral for the involved $k$-forms.
\end{proof}

\subsection{Second step \label{sec:consminsec}}
We now consider the situation of a finite element system $A$ included in a compatible finite element system $B$, such that the following sequences are exact:
\begin{equation}
0 \to A^0_0(T) \to A^1_0(T) \to \cdots \to  A^{\dim T}(T) \to \bbR \to 0.
\end{equation}

For each cell $T$ we will augment the spaces $A^\bs (T)$ in such a way that the extension property holds, while $A^\bs_0(T)$ is preserved (i.e. no trace-free elements will be added). We start by carrying out this process on 1D cells, then on 2D cells, and so forth. We will use the following result, where $\tr$ denotes the trace (pull-back) operator $B^\bullet(T) \to B^\bullet(\partial T)$.

\begin{proposition}\label{prop:mixext}
Let $ l\geq 0$ and suppose $T$ is a cell of dimension $l+1$. Fix $k\leq l$. Any $u \in B^k(\partial T)$ such that  $\rmd u \in \tr A^{k+1}(T)$ has a unique extension $\tilde u\in B^k(T)$ satisfying:
\begin{equation}\label{eq:mixharm}
 \rmd \tilde u \in A^{k+1}(T), \ \rmd \tilde u \perp \rmd A^k_0(T) \myand \tilde u \perp \rmd B^{k-1}_0(T).
\end{equation}
\end{proposition}
\begin{proof}
\emph{Existence.} Suppose first $k<l$. Let $v$ be the $A$-harmonic extension of $\rmd u$. Since $\rmd v$ is $A$-harmonic and zero on the boundary we have $\rmd v = 0$. Choose $u'\in B^k(T)$ such that $\rmd u' = v$. We have $\rmd(\tr u' -u)= 0$.
If $k=0$, $\tr u' - u$ is a constant $c$ and  $u'-c$ is an extension of $u$ with exterior derivative $v$. If on the other hand $k>0$ choose $u''\in B^{k-1}(T)$ such that  $\rmd \tr u'' = \tr u' -u$. Then $u' -\rmd u''$ is an extension of  $u$ with exterior derivative $v$. Adding an element of $\rmd B^{k-1}_0(T)$ we can ensure that it becomes orthogonal to $\rmd B^{k-1}_0(T)$. 

Suppose now $k=l$. Pick $v \in A^{k+1}(T)$ such that $\int v = \int u$ and $v \perp \rmd A^k_0(T)$. Pick $u'\in B^k(T)$ such that $\rmd u' = v$. By Stokes we have $\int(\tr u' - u) = 0$, so we can choose $u''\in B^{k-1}(T)$ such that $\rmd \tr u'' = \tr u' -u$. Then $u' -\rmd u''$ is, as before, an extension of $u$ with exterior derivative $v$ and by adding an element of $\rmd B^{k-1}_0(T)$ we can ensure orthogonality to $\rmd B^{k-1}_0(T)$.

\emph{Uniqueness.} If $\tilde u\in B^k_0(T)$ satisfies (\ref{eq:mixharm}) we have $\rmd \tilde u = 0$. The condition $\tilde u \perp \rmd B^{k-1}_0(T)$ then gives $\tilde u=0$.

\end{proof}

Suppose we have realized our plan for cells of dimension at most $l$ ($l\geq 0$), the augmented FES being denoted $\tilde A$. We consider a cell $T$ of dimension $l+1$.

We equip $\tilde A^k(\partial T)$ with a scalar product, for instance the sum of the local ones, defined for $T' \subcell \partial T$. For $k \leq l$ we put:
\begin{equation}
F^k= \{u \in \tilde A^k(\partial T) \ : \ u \perp \tr A^k(T) \myand \rmd u \in \tr A^{k+1}(T) \},
\end{equation}
and denote by $\tilde F^k$ the extensions of elements of $F^k$ defined by Proposition \ref{prop:mixext} above. We also denote:
\begin{equation}
G^k= \{u \in \tilde A^k(\partial T) \ : \ \rmd u \perp \tr A^{k+1}(T) \myand u \perp \{ v\in  \tilde A^k(\partial T) \ : \ \rmd v = 0 \} \},
\end{equation}
and denote by $\tilde G^k$ the $B$-harmonic extensions of elements of $G^k$. Finally we put:
\begin{equation}\label{eq:tilde}
\tilde A^k(T) = A^k(T) + \tilde F^k + \tilde G^k.
\end{equation}
For $k= l+1$ we take simply $\tilde A^k(T) = A^k(T)$.

\begin{proposition}
On the $(l+1)$-skeleton of $\calT$, $\tilde A$ is a compatible finite element system intermediate between $A$ and $B$ such that  $\tilde A^k_0(T) = A^k_0(T)$ for all $k$ and $T$.
\end{proposition}

\begin{proof}
We remark that $\rmd \tilde F^k \subseteq A^{k+1}(T)$. If $k<l$, $\rmd \tilde G^k \subseteq \tilde F^{k+1}(T)$ and if $k=l$, $\tilde G^k =0$. Therefore $\rmd \tilde A^k(T) \subseteq \tilde A^{k+1}(T)$. We have a direct sum decomposition:
\begin{equation}
\tilde A^k(\partial T) = \tr A^k(T) \oplus F^k \oplus G^k.
\end{equation}
It follows in particular that $\tr : \tilde A^k(T) \to \tilde A^k(\partial T)$ is surjective. Moreover if $ u \in A^k(T)$, $ v\in \tilde F^k $ and $w \in \tilde G^k$ satisfy $\tr (u + v + w) = 0$, then $\tr v =0$ and $\tr w =0$, which gives $v=0$ and $w=0$. If in fact  $u + v + w = 0$, we will also have $u=0$. Therefore the sum (\ref{eq:tilde}) is direct and:
\begin{equation}
\tilde A^k_0(T) =  A^k_0(T).
\end{equation}
This completes the proof.
\end{proof}

The construction is therefore complete, and provides a minimal compatible finite element system containing $A$, by Corollary \ref{cor:minnec}.

\section{Minimality of some finite element systems \label{sec:ex}}

We now consider more concrete examples, based on polynomial differential forms on simplices or hypercubes. 
The two steps of our general construction enter differently in them:
\begin{itemize}
\item In Section \ref{sec:mintrim}, on trimmed differential forms, the local cohomology groups $ \rmH^k(A^\bs_0(T))$ are non-trivial, but once this problem is fixed, the obtained spaces turn out to have the extension property.
\item In Section \ref{sec:minser}, on serendipity elements,  the local cohomology groups are trivial, so one only needs to add spaces of forms that ensure the extension property.
\item In Section \ref{sec:mintnt}, on TNT elements, the extension property holds for the spaces one starts with, but local sequence exactness fails.  By adding tiny bubble functions $E^k(T)$ to fix the cohomology requirement, one also needs to add forms that extend them to higher-dimensional cells.
\item In Section \ref{sec:minsp}, in which we introduce ``small pleasures'' elements, the sequence one starts with has nontrivial cohomology and does not satisfy the extension property.
\end{itemize} 

\subsection{Minimality of trimmed polynomial differential forms\label{sec:mintrim}}

We now adopt notations from finite element exterior calculus~\cite{ArnFalWin06}. 
We first consider finite element systems on a simplicial complex $\calT$. 
Let $\poly_r \Lambda^k(T)$ denote the space of $k$-forms on $T\in\calT$ whose coefficients are polynomials of degree at most $r$. The space $\poly_r^-\Lambda^k(T)$ is defined with the help of the operator $\kappa$ on forms, which is contraction by the vector field $x \mapsto x$ (different choices of origin yield the same space). Then:
\begin{equation}
\poly_r^-\Lambda^k(T) = \{ u \in \poly_r \Lambda^k(T) \ : \  \kappa u \in  \poly_r \Lambda^{k-1}(T) \}.
\end{equation}

In keeping with our previous notations, $\poly_r^-\Lambda^k_0(T)$ is the subspace of $\poly_r^-\Lambda^k(T)$ consisting of differential forms whose pullback to the boundary $\partial T$ of $T$ is $0$. We also use freely that these spaces form exact sequences under the exterior derivative, which follows from the fact that they constitute a \emph{compatible} FES. The underlying facts are proved in standard finite element language in \cite{ArnFalWin06}, but it is also checked with the general tools of finite element systems in \cite{ChrRap15}.

Similar statements hold for the $\poly_r\Lambda^k(T)$ spaces, but here sequence exactness requires decreasing $r$ through the complex. That is, under the exterior derivative, the sequence:
\begin{equation}
\poly_{r+1} \Lambda^{k-1}_0(T) \to \poly_r \Lambda^k_0(T) \to \poly_{r-1} \Lambda^{k+1}_0(T), 
\end{equation}
is exact, whereas the sequence:
\begin{equation}
\poly_{r} \Lambda^{k-1}_0(T) \to \poly_r \Lambda^k_0(T) \to \poly_{r} \Lambda^{k+1}_0(T), 
\end{equation}
in general is not. The cohomology group of the latter sequence is denoted:
\begin{equation}
\rmH^{k} \poly_{r}\Lambda^\bs_0(T).
\end{equation}

\begin{lemma} We have :
\begin{equation}
\rmd \poly_r^-\Lambda^k_0(T) = \rmd \poly_r\Lambda^k_0(T) .
\end{equation}
\end{lemma}
\begin{proof}
Pick $u \in \poly_r\Lambda^k_0(T)$ and set $v:= \rmd u$.
Then $v \in \poly_{r-1}\Lambda^{k+1}_0(T)\subseteq \poly_{r}^-\Lambda^{k+1}_0(T)$ and $\rmd v = 0$. Hence there is $w \in \poly_{r}^-\Lambda^{k}_0(T)$ such that $\rmd w = v$.
\end{proof}

The following identity shows that  the spaces $\poly_r^-\Lambda^k(T) $ constitute a \emph{minimal} compatible FES containing the spaces $\poly_{r-1}\Lambda^k(T)$.
\begin{proposition} We have:
\begin{equation}
\dim  \poly_r^-\Lambda^k_0(T) = \dim  \poly_{r-1}\Lambda^k_0(T) + \dim \rmH^{k+1} \poly_{r-1}\Lambda^\bs_0(T).
\end{equation}
\end{proposition}
\begin{proof}
We use the vertical line notation $|$ to denote restriction of an operator to a subspace.
We have:
\begin{align}
\dim  \poly_r^-\Lambda^k_0(T) &= \dim \rmd \poly_r^-\Lambda^k_0(T) + \dim \ker \rmd | \poly_r^-\Lambda^k_0(T),\\
& = \dim \rmd \poly_r\Lambda^k_0(T) + \dim \rmd \poly_r^-\Lambda^{k-1}_0(T),\\
&  =  \dim \rmd \poly_r\Lambda^k_0(T)  + \dim \rmd \poly_r\Lambda^{k-1}_0(T),
\end{align}
and:
\begin{align}
\dim  \poly_{r-1}\Lambda^k_0(T) &= \dim \rmd \poly_{r-1}\Lambda^k_0(T) + \dim \ker \rmd | \poly_{r-1} \Lambda^{k}_0(T),\\
&= \dim \rmd \poly_{r-1}\Lambda^k_0(T) + \dim \rmd \poly_r\Lambda^{k-1}_0(T).
\end{align}
Moreover:
\begin{align}
\dim \rmH^{k+1} \poly_{r-1}\Lambda^\bs_0(T) & = \dim \ker \rmd | \poly_{r-1}\Lambda^{k+1}_0(T)  - \dim \rmd \poly_{r-1}\Lambda^k_0(T),\\
&= \dim \rmd \poly_r\Lambda^k_0(T)  - \dim  \rmd \poly_{r-1}\Lambda^k_0(T).
\end{align}
The result follows by addition of the above equalities.
\end{proof}

\subsection{Minimality of serendipity elements\label{sec:minser}}

In \cite{ArnAwa14}, Arnold and Awanou define spaces of \textit{serendipity} finite element differential forms on $n$-dimensional cubes, denoted $\calS_r\alt^k(I^n)$.
They prove a subcomplex property (Theorem 3.3) and trace property (Theorem 3.5), which is equivalent, in our notation, to showing that the spaces $\calS_{r-k} \alt^k(I^n)$
constitute a FES, for any fixed $r\geq n$. Arnold and Awanou also define these spaces for $r < n$, but then one does not get a full complex, but one that stops at $k = r$. In many applications this is useful, since one does not require the full complex but rather two consecutive spaces, possibly three. But the FES framework has not been designed for this.

They also prove that the spaces $\poly_{r -k - 2(n-k)} \alt^{n-k} (I^n)$ provide unisolvent degrees of freedom (Theorem 3.6).  By \cite[Proposition 2.5]{ChrRap15}, this guarantees that the extension property holds. They also check the sequence exactness property (page 1566). Thus, in our language, they have shown that the serendipity finite element spaces constitute a compatible finite element system. The inclusion $\poly_{r-k}\alt^k(I^n) \subseteq A^k(I^n)$ is obvious from their construction and it is implicit from their Proposition 3.7 that $A^k_0(I^n) = \poly_{r-k}\alt^k_0(I^n)$.

Thus one can interpret many of their results, by the statement that the serendipity elements define a minimal compatible finite element system containing $ \poly_{r-k}\alt^k(I^n)$. We now explain this, proving the main statements above, in our framework.

It should be noted that the construction in \cite{ArnAwa14} is very different from the one we propose in Section \ref{sec:consminsec} (the first step in our construction is not necessary here, for reasons provided below). Our construction depends on a choice of compatible finite element system $B$ containing $\poly_{r-k}\alt^k(I^n)$. For that purpose we could take the spaces:
\begin{equation}
B^\bs(I^n) = \poly^-_{r} \alt^\bs(I)  \otimes \ldots \otimes \poly^-_{r} \alt^\bs(I).
\end{equation}
We would also need a scalar product. One can use $\rmL^2$ scalar products, or more algebraic ones, for instance one that makes some preferred basis in $B$ orthonormal.

\begin{remark}
We will take $I$ to be the interval $[0, 1]$ while Arnold and Awanou take it to be $[-1,1]$.
Also, we will only address the cases $r\geq n$, although the arguments can be adapted to treat lower order cases as well.
\end{remark}

\begin{lemma}\label{lem:zerotens}
Let $B_\alpha^k$ be a space of differential $k$-forms on $I^{n-1}$, for each $0 \leq \alpha \leq r$ and for each $k$.
Define the sum of graded tensor product spaces of forms on $I^n$:
\begin{equation}
A^\bs = \sum_{\alpha = 0}^r B_\alpha^\bs \otimes \poly_{\alpha} \alt^{\bs}(I).
\end{equation}
Consider $\delta_n^0$ and $\delta_n^1$, the $n$-th trace operators, defined respectively as pullbacks by the injections $ I^{n-1} \times \{0\} \to I^n$   and $I^{n-1} \times  \{1\} \to I^n$. Then the intersections of the kernels of $\delta_n^0$ and $\delta_n^1$ on $A^\bs$ is :
\begin{equation}
(\ker \delta_n^0 | A^\bs ) \cap (\ker \delta_n^1 | A^\bs)  = \sum_{\alpha = 0}^r B_\alpha^\bs \otimes \poly_{\alpha} \alt_0^{\bs}(I).
\end{equation}
\end{lemma}
\begin{proof}
We let $\lambda_0$ (resp. $\lambda_1$) denote the element of $\poly_{1} \alt^{0}(I)$ taking the values $1$ (resp. $0$) at $0$ and $0$ (resp. $1$) at $1$. Observe that for $\alpha \geq 1$ we have a direct sum decomposition:
\begin{equation}
\poly_{\alpha} \alt^{\bs}(I) = \bbR \lambda_0 \oplus \bbR \lambda_1 \oplus \poly_{\alpha} \alt^{\bs}_0(I).
\end{equation}
We also have:
\begin{equation}
\poly_{0} \alt^{\bs}(I) = \bbR (\lambda_0  + \lambda_1) +  \bbR \rmd \lambda_0.
\end{equation}

Now take an element $u$ of $A^k$ and write it as a sum:
\begin{equation}
u =     v_0 \otimes  (\lambda_0  + \lambda_1) + v'_0 \otimes \rmd \lambda_0 + \sum_{\alpha = 1}^r (v_{\alpha_0} \otimes \lambda_0 + v_{\alpha_1} \otimes \lambda_1 + v_\alpha),
\end{equation}
with:
\begin{equation}
v_0  \in B_0^k,\ v'_0  \in B_1^{k-1}, \textrm{ and } v_\alpha  \in \sum_l B_\alpha^{k-l} \otimes \poly_{\alpha} \alt^{l}_0(I).
\end{equation}
Observe that if $u$ is in the intersection of the kernels, we have:
\begin{align}
0 = \delta_n^0 u & = v_0 + \sum_{\alpha = 1}^r v_{\alpha_0},  \\
0 = \delta_n^1 u & = v_0 + \sum_{\alpha = 1}^r v_{\alpha_1}.
\end{align}
Then we are left with:
\begin{equation}
u = v'_0 \otimes \rmd \lambda_0  + \sum_{\alpha = 0}^r  v_\alpha,
\end{equation}
and this completes the proof.
\end{proof}

\begin{proposition}\label{prop:polyzc}
On $I^n$ let $x_j$ denote the $j$-th coordinate function. We have:
\begin{equation}\label{eq:polyrz}
\poly_r \alt^k_0(I^n) =  \lspan \left\{ u \prod_{j \not \in J} x_j (1- x_j) \rmd x_J \ :  \begin{array}{l} u \in \poly_{r - 2(n-k)}(I^n), \\  J \subseteq \{1, \ldots, n\}, \\  \card J = k. \end{array}\right\},
\end{equation}
where, for $J =\{ j_1, \ldots, j_k \}$ we define $\rmd x_J = \rmd x_{j_1} \wedge \ldots \wedge \rmd x_{j_k}$.
\end{proposition}
\begin{proof}
We have:
\begin{equation}
\poly_r \alt^\bs(I^n) = \sum_{\alpha_1 + \ldots + \alpha_n = r} \poly_{\alpha_1}\alt^\bs (I) \otimes \ldots \otimes \poly_{\alpha_n}\alt^\bs (I), 
\end{equation}
where we sum over integer $\alpha_i$ in the range $[0, r]$.

In order to determine $\poly_r \alt^\bs_0(I^n)$ we apply Lemma \ref{lem:zerotens}, or its equivalent, in each coordinate direction. We get:
\begin{equation}
\poly_r \alt^\bs_0(I^n) = \sum_{\alpha_1 + \ldots + \alpha_n = r} \poly_{\alpha_1}\alt^\bs_0(I) \otimes \ldots \otimes \poly_{\alpha_n}\alt^\bs_0(I), 
\end{equation}
Now, concerning each factor, we have with $x$ the coordinate function on $I$:
\begin{align}
\poly_{q}\alt^0_0(I) & = \{ u(x) x(1-x) \ : \ u \in \poly_{q-2}(I)  \},\\
\poly_{q}\alt^1_0(I)& = \poly_{q}\alt^1(I).
\end{align}

It follows that the left hand side of (\ref{eq:polyrz}) is included in the right hand side. The inclusion in the other direction is trivial, so the proof is complete.
\end{proof}

\begin{remark}
The result just shown seems to be taken as self-evident in \cite[page 1569]{ArnAwa14}.
\end{remark}

\begin{proposition}\label{prop:zerodual}
The integrated wedge product defines an invertible bilinear form:
\begin{equation}
(\int \cdot \wedge \cdot) \ : \  \poly_r \alt^k_0(I^n) \times \poly_{r - 2(n-k)} \alt^{n-k} (I^n) \to \bbR.
\end{equation}
\end{proposition}
\begin{proof}
It follows from the preceding proposition that the two spaces have the same dimension.

Fix $u \in \poly_r \alt^k_0(I^n)$. We want to show that the map sending $v\mapsto \displaystyle \int u\wedge v\in\bbR$ is non-zero only if $u$ is zero. Write $u$ in the form:
\begin{equation}
u = \sum_J \left(u_J \prod_{j \not \in J} x_j (1- x_j) \rmd x_J\right),
\end{equation}
with:
\begin{equation}
u_J \in \poly_{r - 2(n-k)}(I^n),\ J \subseteq \{1, \ldots, n\} \textrm{ and } \card J = k.
\end{equation}
Letting $J'$ denote the complement of $J$ in $\{1, \ldots, n\}$ we define:
\begin{equation}
v= \sum_J \epsilon_J u_J \rmd x_{J'} \in \poly_{r - 2(n-k)} \alt^{n-k} (I^n).
\end{equation}
Then we have:
\begin{equation}
u \wedge v = \sum_J \epsilon_J u_J^2 \prod_{j \not \in J} x_j (1- x_j)\rmd x_J \wedge \rmd x_{J'} .
\end{equation}
We then choose the signs $\epsilon_J$ such that $\epsilon_J \rmd x_J \wedge \rmd x_{J'} $ is the volume form on $I^n$. Then $\int u\wedge v$ is zero iff $u$ is zero.
\end{proof}

\begin{proposition}\label{prop:zeroce} Fix $r \geq 0$.

For $0 \leq k < n$, the following sequence, for the exterior derivative, is exact:
\begin{equation}
\poly_{r+1} \alt^{k-1}_0(I^n) \to \poly_r \alt^k_0(I^n) \to \poly_{r-1} \alt^{k+1}_0(I^n).
\end{equation}

The sequence:
\begin{equation}
\poly_{r+1} \alt^{n-1}_0(I^n) \to \poly_r \alt^n_0(I^n) \to \bbR,
\end{equation}
is also exact (where the last arrow is integration).

\end{proposition}
\begin{proof}
Consider the first statement. By Proposition \ref{prop:zerodual}, the sequence is paired with the following one:
\begin{equation}
\poly_{r-1- 2(n-k)} \alt^{n-k+1}(I^n) \from \poly_{r - 2(n-k)} \alt^{n-k}(I^n) \from \poly_{r+1- 2(n-k)} \alt^{n-k-1}(I^n).
\end{equation}
The differential in this sequence is again the exterior derivative. Recall also that:
\begin{equation}
\int_T \rmd u \wedge v = \pm \int_T u \wedge \rmd v,
\end{equation}
when (for instance) $u$ is zero on $\partial T$.
Since the latter sequence is exact, the former one must be too.

Consider now the second statement. The sequence is paired with:
\begin{equation}
\poly_{r-1} \alt^{n-1}(I^n) \from \poly_{r} \alt^{0}(I^n) \from \bbR.
\end{equation}
On the right we have inclusion of the constants. This sequence is exact.
\end{proof}

We may summarize our results as follows:
\begin{proposition}
Fix $r \geq n$. Let $A^k(I^m)$ be a compatible finite element system on hypercubes of dimension at most $n$, such that:
\begin{equation}\label{eq:sercp}
\poly_{r-k}\alt^k(I^m) \subseteq A^k(I^m).
\end{equation}
Then $A^k(I^m)$ is a minimal compatible finite element system with this property if and only if:
\begin{equation}\label{eq:minpz}
A^k_0(I^m) = \poly_{r-k}\alt^k_0(I^m). 
\end{equation}
In this case, the spaces $\poly_{r +k - 2m} \alt^{m-k} (I^m)$ provide unisolvent degrees of freedom on $A$.
\end{proposition}

\begin{proof}
The first statement follows from the dimensional characterization of minimality provided by Corollary~\ref{cor:minnec}, given the exactness property proved in Proposition \ref{prop:zeroce}. The second statement recalls Proposition~\ref{prop:zerodual}.
\end{proof}

Again we point out that \cite{ArnAwa14} and Section \ref{sec:con} use quite different means to define a compatible FES $A$ such that (\ref{eq:sercp}) and (\ref{eq:minpz}) hold.

\subsection{Minimality of TNT elements\label{sec:mintnt}}

We now consider the work of Cockburn and Qiu regarding the `TiNiest spaces containing Tensor product spaces of polynomials' on cubes or \emph{TNT} elements for short~\cite{CocQiu14}.
We will show how to derive their minimality claims from the framework of FES. Moreover we cast their constructions in the language of differential forms, which is uniform in the treatment of dimensions and degree of forms. We will also make explicit use of other facts from homological algebra such as the Kunneth theorem.

Fix $r\geq 1$.  Define the spaces, on hypercubes:
\begin{equation}
\label{eq:a-tnt-def}
A^\bs(I^n) = \calQ_r\alt^\bs(I^n) = \poly_r \alt^\bs(I)  \otimes \ldots \otimes \poly_r \alt^\bs(I).
\end{equation}
We remark that on an interval $\poly_r \alt^\bs(I)$ has the extension property. Therefore, by Proposition \ref{prop:tensorfes}, $A$ also has the extension property. However the sequence exactness fails. We want to define a minimal compatible finite element system containing $A$. This is accomplished for dimension $n=3$ in \cite{CocQiu14}. 
 
By \cite[Lemma 3.10]{Chr09AWM} we have:
\begin{equation}
A^\bs_0(I^n) = \poly_r \alt^\bs_0(I)  \otimes \ldots \otimes \poly_r \alt^\bs_0(I).
\end{equation}

Looking at the sequence:
\begin{equation}\label{eq:zerocohr}
0 \to \poly_r \alt^0_0(I) \to \poly_r \alt^1_0(I) \to 0,
\end{equation}
we notice that there is a non-zero cohomology group only in the 1-forms, and that it has dimension $2$. We write out the Kunneth theorem:
\begin{equation}
\rmH^\bs A^\bs_0(I^n) \approx \rmH^\bs (\poly_r \alt^\bs_0(I)) \otimes \ldots \otimes \rmH^\bs (\poly_r \alt^\bs_0(I)).
\end{equation}
Since, as we remarked above, $\rmH^0 (\poly_r \alt^\bs_0(I)) = 0$, most of the terms cancel. There is only one non-zero cohomology group in $A^\bs(I^n)$, namely the one in $n$-forms. Explicitly:
\begin{align}
\rmH^{k} A^\bs_0(I^n) & = 0,\ \textrm{ for } k < n, \label{eq:aexn}\\
\rmH^n A^\bs_0(I^n) & \approx \rmH^1 (\poly_r \alt^\bs_0(I)) \otimes \ldots \otimes \rmH^1 (\poly_r \alt^\bs_0(I)), \label{eq:hnaz}
\end{align}
and the latter space has dimension $2^n$.

Now, the sequence that we want to make exact is really:
\begin{equation}
0 \to A^0_0(I^n) \to A^1_0(I^n) \to \ldots \to A^n(I^n) \to \bbR \to 0,
\end{equation}
where the second to last arrow is integration. Thus the notation in (\ref{eq:hnaz}) is slightly different from the one used in the beginning of the article, e.g in (\ref{eq:mindimeq}), the difference being whether we consider that there is an $\bbR$ in the end  of $A^\bs_0(I^n)$ or not. We notice that $A^n(I^n)$ contains elements with non-zero integral. We thus have:

\begin{proposition}\label{prop:cocmincond}
Let $B$ denote a compatible finite element system containing $A$, as defined in (\ref{eq:a-tnt-def}). Then $B$ is a minimal one containing $A$ if and only if, for all $k$ and $n$:
\begin{align}
B^k_0(I^n) & = A^k_0(I^n), \ \textrm{ for } k \neq  n-1 ,\\
\dim B^{n-1}_0 (I^n) & = \dim A^{n-1}_0(I^n)  + 2^n -1.
\end{align}
\end{proposition}
This provides the dimension count for any minimal finite element system containing $A$, using Proposition \ref{prop:extdim}.

Any construction of a minimal compatible finite element system $B$ containing $A$ would need to add spaces $E^{n-1}(I^n)$ to $A^{n-1}(I^n)$ such that:
\begin{align}
\dim E^{n-1}(I^n) & = 2^n -1,\\
B^{n-1}_0(I^n) & = A^{n-1}_0(I^n) \oplus E^{n-1}(I^n),
\end{align}
as in the `first step' outlined in Section \ref{sec:consminfir}. The second step, outlined in Section \ref{sec:consminsec},  would have to ensure that the extension property holds, after these additions, since they also take place on the faces of $I^n$.

To be more specific, we notice that if $x$ denotes the coordinate function on the interval $I$, we can choose two generators of $\rmH^1 (\poly_r \alt^\bs_0(I))$, to be $\rmd x $ and $P(x)\rmd x$, where $P$ is a polynomial of degree $r$ with integral $0$ on $I$. We also impose $P(1-x) = P(x)$ if $r$ is even and $ P(1-x) = - P(x)$ if $r$ is odd (to ensure that our construction becomes invariant under reflections in the coordinate directions). Then we let $Q$ be the polynomial of degree $r+1$ such that $Q(0)= Q(1) = 0$ and $Q' = P$.

Denote by $x_i$ the $i$-th coordinate map on $I^n$. Let $J$ denote a non-empty subset of $\{1, \ldots, n\}$. Define the $n$-form $f_J$ on $I^n$:
\begin{equation}
f_J = \bigotimes_{j \in J} P(x_j) \rmd x_j \bigotimes_{j \not \in J} \rmd x_j,
\end{equation}
where the notation is taken to mean that we have the form $P(x) \rmd x$ in the directions of $J$ and the form $\rmd x$ in the other ones (thus there is a permutation in the indices going on). Since $J$ is non-empty, $f_J$ has integral $0$. Notice that we have thus defined $2^n-1$ linearly independent forms in $A^n(I^n)$ with zero integral.

We now observe that:
\begin{equation}
f_J = \frac{1}{\card J} \rmd g_J, \textrm{ with } g_J = \sum_{i \in J}  Q(x_i)\bigotimes_{j \in J\setminus \{i\}} P(x_j) \rmd x_j \bigotimes_{j \not \in J} \rmd x_j.
\end{equation}
This notation is taken to mean that we can obtain $g_J$ from $f_J$ by replacing, for each index $i$ in $J$, the term $P(x_i)\rmd x_i$ by $Q(x_i)$, and summing over $i \in J$.

We define:
\begin{equation}
E^{n-1}(I^n) = \lspan \{ g_J \ : \  J \subseteq \{1, \ldots, n\}, \ J \neq \emptyset \}.
\end{equation}
We remark that:
\begin{lemma}\label{lem:egood}
The elements of $E^{n-1}(I^n)$ are zero on $\partial I^n$ and the map defined in  (\ref{eq:disoe}) is an isomorphism.
\end{lemma}

Thus they are adequate for the first step of our general construction. Notice also that  these spaces are invariant under permutation of indices and reflections $x_i \mapsto 1-x_i$.

For the second step, one should keep in mind that we have added spaces $E^{d-1}(T)$ for all faces $T$ of dimension $d$ in $I^n$, and that the extension property is therefore no longer guaranteed. We could proceed by our general method. But inspired by \cite{CocQiu14}, we prefer the following explicit approach:

We determine a face $T$ of $I^n$ by a choice of indices $J_T\subseteq \{1, \ldots , n \}$ and for the indices $i \in \{1, \ldots , n \}\setminus J_T$ a number $N_T(i) \in \{0, 1 \}$. Then the associated face $T$ is:
\begin{equation}
T = \{ (x_1, \ldots , x_n ) \  : \ \textrm{ for }  i \not \in J_T \quad x_i = N_{T}(i) \}.
\end{equation} 
Let $T$ be a $d$ dimensional face of $I^n$. An element of $E^{d-1} (T)$ is extended to $I^n$ as follows. Let $u$ be the element in question. Its extension is chosen to be:
\begin{equation}\label{eq:tildeu}
\tilde u = u \bigotimes_{i \not \in J_T} x_i^{N_T(i)}(1-x_i)^{1-N_T(i)}.
\end{equation}
We then define:
\begin{equation}\label{eq:defcoc}
B^k(I^n) = A^k(I^n) + \sum_{\dim T = k+1} \tilde E^k(T),
\end{equation}
where the sum is over $k+1$ dimensional faces of $I^n$ and $\tilde E^k(T)$ denotes the space of extensions of elements of $E^k(T)$ defined by $(\ref{eq:tildeu})$.

\begin{proposition}
The definition (\ref{eq:defcoc}) provides a minimal compatible FES containing $A$.
\end{proposition}
\begin{proof}
Suppose that (\ref{eq:tildeu}) holds for some $u \in E^k(T)$.  Then:
\begin{align}
\rmd \tilde u   & =~ u \bigotimes_{i \not \in J_T} x_i^{N_T(i)}(1-x_i)^{1-N_T(i)}~ +\notag \\
 & \qquad\quad \sum_{j \not \in J_T} u \bigotimes_{i \not \in J_T}  x_i^{N_T(i)}(1-x_i)^{1-N_T(i)} (\pm \rmd x_i)_{i= j},\\
& \in ~~A^{k+1}(I^n)  + \sum_{\dim T = k+2} \tilde E^{k+1}(T).
\end{align}
Thus $B$ satisfies the requirements to be an element system.
Let $T$ be a $k+1$ dimensional face of $I^n$ and let $T'$ be another one. Let's take an element $u$ of $ E^k(T)$ and restrict $\tilde u$ to $T'$. Choose $i \in J_{T} \setminus J_{T'}$.  The restrictions of $u$ to the faces $x_i = 0$ and $x_i = 1$ inside $T$ are $0$. Therefore the restriction of $\tilde u$ to $T'$ is $0$.
From this we deduce that the sum in (\ref{eq:defcoc}) is direct, and that:
\begin{align}
B^k_0(I^n) & = A^k_0(I^n), \ \textrm{ for } k \neq  n-1 ,\\
B^{n-1}_0 (I^n) & = A^{n-1}_0(I^n)  \oplus E^{n-1}(T) .
\end{align}
Recalling that $A$ has the extension property, we see in particular that:
\begin{align}
\dim B^k(I^n) & = \dim A^k(I^n) +  \sum_{\dim T = k+1} \dim E^k(T),\\
& = \sum_{\dim T \neq k+1} A^k_0(T) + \sum_{\dim T = k+1} \left(\dim A^k_0(T) + \dim E^k(T)\right) ,\\
& = \sum_{T} \dim B^k_0(T), 
\end{align}
where $T$ runs through the set of faces of $I^n$. This proves that $B$ has the extension property by Proposition \ref{prop:extdim}.
The sequence exactness in the form (\ref{eq:coh0}) has already been proved, see (\ref{eq:aexn}) and Lemma \ref{lem:egood}.
Thus we have proved that $B$ is a compatible finite element system.
It is minimal containing $A$ by application of Proposition \ref{prop:cocmincond}.

\end{proof}

If we want degrees of freedom for this compatible FES that yield commuting interpolators, we can take the ones defined in \cite[Proposition 2.8]{ChrRap15}.

\subsection{Dimension count for ``small pleasures'' elements\label{sec:minsp}}

Neither the construction of Arnold and Awanou nor the one of Cockburn and Qiu resolves the question that seems most natural to us, namely, defining the minimal compatible FES on hypercubes containing the sequences
\begin{equation}
\label{eq:minseqPrLkIn}
\poly_r\alt^0(I^n) \to \poly_r\alt^1(I^n) \to \ldots \to \poly_r\alt^n(I^n).
 \end{equation}
The preceding general construction provides such a mcFES. Being neither the serendipity elements, nor the tiniest elements in the Cockburn-Qiu sense, we refer to these as the ``Small Pleasures'' elements.
We provide a dimension count of the Small Pleasures elements with the following result.

\begin{lemma}
\begin{equation}
\label{eq:dimHkPrL0In}
\dim \rmH^k \poly_r \alt^\bs_0(I^n)  = { r + 2k-n -1 \choose k-1}{r +k -n -1 \choose n - k}.
\end{equation}
\end{lemma}
\begin{proof}

We have an exact sequence:
\begin{align}
0 \to \poly_{r+k-1} \alt^0_0(I^n) & \to  \poly_{r+k-2}\alt^1_0(I^n) \to \ldots \notag \\
& \ldots \to \poly_{r+1} \alt^{k-2}_0(I^n)   \to \poly_r\alt^{k-1}_0(I^n) \to \rmd \poly_r\alt^{k-1}_0(I^n) \to 0.
 \end{align}
This gives:
\begin{equation}
\dim \rmd \poly_r\alt^{k-1}_0(I^n) = \sum_{i= 0}^{k-1} (-1)^i \dim \poly_{r+i}\alt^{k-1 -i }_0(I^n).
\end{equation}

We also have:
\begin{equation}
\ker \rmd  | \poly_r \alt^k_0(I^n)  =  \rmd  \poly_{r+1} \alt^{k-1}_0(I^n).
\end{equation}
Hence, by the same formula as previously:
\begin{equation}
\dim \ker \rmd | \poly_r\alt^{k}_0(I^n) = \sum_{i= 0}^{k-1} (-1)^i \dim \poly_{r+1+i}\alt^{k-1 -i }_0(I^n).
\end{equation}

Therefore:
\begin{align}
\dim \rmH^k \poly_r \alt^\bs_0(I^n) & = \sum_{i= 0}^{k-1} (-1)^i (\dim \poly_{r+1+i}\alt^{k-1 -i }_0(I^n) - \dim \poly_{r+i}\alt^{k-1 -i }_0(I^n)),\\
& = \sum_{i= 0}^{k-1} (-1)^i (\dim \poly_{r+1+i  - 2(n-k+1+i)}  \alt^{n -k+ 1 + i }(I^n) -\\
& \qquad \dim \poly_{r+i - 2(n-k+1+i)}\alt^{n-k+1 +i }(I^n)),\\
& = \sum_{i=0}^{k-1} (-1)^i \dim \calH_{r-1-i  - 2(n-k)}  \alt^{n -k+ 1 + i }(I^n). 
\end{align}

Now, following \cite[Section 3.2]{ArnFalWin06}, the Koszul sequence on $I^n$:
\begin{equation}
0 \to \calH_t \alt^n \to \calH_{t+1} \alt^{n-1} \to \ldots \to \calH_{t+ k-1} \alt^{n-k+1} \to \kappa \calH_{t+ k-1} \alt^{n-k+1} \to 0,
\end{equation}
is exact. We take $t= r+ k -2n$, and deduce:
\begin{align}
 \sum_{i = 0}^{k-1} (-1)^i \dim\calH_{t+k-1-i}\alt^{n-k+1+i} & =  \dim \kappa \calH_{t+ k-1} \alt^{n-k+1},\\
& = { n + t + k -1 \choose k - 1} {t +n -1 \choose n-k},\\
& = { r + 2k-n -1 \choose k-1}{r +k -n -1 \choose n - k}.
\end{align}

This concludes the proof.
\end{proof}

Using (\ref{eq:dimHkPrL0In}) and (\ref{eq:mindimeq}), we can compute the dimension of $B^k_0(I^n)$ where $B$ is a minimal compatible finite element system containing (\ref{eq:minseqPrLkIn}).
We give some some sample values in Table~\ref{tab:dimcount}.
These spaces will be investigated further in future work.

\begin{table}
\[
\begin{tabular}{l|ccccccc}
 & $r=4$ & $r=5$ & $r=6$ & $r=7$ & $r=8$ & $r=9$ & $r=10$ \\
 \hline
$k=0$ & 0 & 1 & 4 & 10 & 20 & 35 & 56 \\
$k=1$ & 11 & 27 & 54 & 95 & 153 & 231 & 332 \\
$k=2$ & 45 & 81 & 133 & 204 & 297 & 415 & 561 \\
$k=3$ & 35 & 56 & 84 & 120 & 165 & 220 & 286
\end{tabular}
\]
\caption{Some computed values for $\dim B^k_0(I^3)$ where $B$ is a minimal compatible finite element system containing (\ref{eq:minseqPrLkIn}).}
\label{tab:dimcount}
\end{table}

\section*{Acknowledgements}
SHC was supported by the European Research Council through the FP7-IDEAS-ERC Starting Grant scheme, project 278011 STUCCOFIELDS.
AG was supported in part by NSF Award 1522289.

%\bibliography{shc-akg}
\bibliography{../Bibliography/alexandria,../Bibliography/newalexandria,../Bibliography/mybibliography}{}
\bibliographystyle{abbrv}

\end{document}